\documentclass[12pt,oneside,final]{amsart}
\usepackage[cp1250]{inputenc}
\usepackage[T1]{fontenc}
\usepackage{amssymb}
\usepackage{color}
\renewcommand{\thefootnote}{}
\newtheorem{twr}{Twierdzenie}
\newtheorem{thr}[twr]{Theorem}
\newtheorem{prp}[twr]{Proposition}

\newtheorem{lm}[twr]{Lemma}

\newtheorem{crl}[twr]{Corollary}
\renewcommand{\thefootnote}{}
\newcommand{\da}{\partial}
\newcommand{\db}{\bar\partial}
\newcommand{\dc}{i{\partial}\bar\partial}
\newcommand{\ta}{\theta}
\newcommand{\tb}{\bar\theta}
\newcommand{\tc}{\theta\bar{\theta}}

\newcommand{\po}{{\partial\Omega}}

\newcommand{\ma}{{(i\partial\bar\partial u)^2}}

\numberwithin{twr}{section}
\title[Monge-Amp\`ere operator]{Monge-Amp\`ere operator on four dimensional almost complex manifolds   }
\author[S. Pliś]{Szymon Pliś}
\address{  Institute of Mathematics, Cracow University of Technology, Warszawska 24, 31-155
    Kraków, Poland
}

\email{splis@pk.edu.pl}
\subjclass[2010]{ 32W20,32Q60, 32U05, 32U40, 35J60}
\keywords{  Monge-Amp\`ere operator, almost complex manifold, plurisubharmonic function, Richberg theorem}

\begin{document}\thispagestyle{empty} \footnotetext{The  author was partially supported by the NCN grant 2011/01/D/ST1/04192.}
\renewcommand{\thefootnote}{\arabic{footnote}} 

\begin{abstract}
We define the Monge-Amp\`ere operator $\ma$ for continuous $J $-plurisubharmonic functions on four dimensional almost complex manifolds.
\end{abstract}

\maketitle
\section{introduction}
Recently several papers concerning  the  plurisubharmonic functions on almost complex manifolds has appeared (see for example \cite{p}, \cite{h-l} and \cite{k}). A very useful tool to work with plurisubharmonic functions on complex manifolds  is the complex Monge-Amp\`ere operator $(\dc u)^n$, which is well defined  for not necessary smooth plurisubharmonic functions  (see \cite{b-t1} and \cite{b-t3}). It also seems  convenient to define this operator on almost complex manifolds. In this paper we deal with that problem in the  basic case of continuous plurisubharmonic functions on four real dimensional almost complex manifolds. The main theorem is the following:
\begin{thr}\label{mine theorem}
Let $M$ be a four dimensional almost complex manifold. Then the Monge-Amp\`ere operator $\ma$ is well defined as a regular Borel measure for any continuous plurisubharmonic function $u$.
\end{thr}
The construction of $\ma$ is given  in  subsection \ref{definition}, where we prove Theorem \ref{MAdlaciaglych}, which is a slight generalisation of Theorem \ref{mine theorem}.

In the integrable case a key property which allows to define a wedge product of plurisubharmonic currents $\dc u$  is their positivity and closeness. N. Pali showed in \cite{p} that in our case such currents are positive but not closed, and it is the main difficulty in the construction of the operator.

An important step to define the Monge-Amp\`ere operator is a result, which is interesting by itself, about the smoothing of continuous plurisubharmonic functions. We prove it (in any dimension) in section~ ~\ref{approximation}.


\section{preliminaries}
\subsection{almost complex manifold}
We say that $(M,J)$ is an almost complex manifold if $M$ is a manifold and $J$ is an ($\mathcal{C}^\infty$ smooth) endomorphism of the tangent bundle $TM$, such that $J^2=-\rm{id}$. The real dimension of $M$ is even in that case. We will always denote by $n$ a complex dimension of $M$: $n=\dim_\mathbb{C}M=\frac{1}{2}\dim_\mathbb{R}M$.

 We have a direct sum decomposition $T_{\mathbb{C}}M=T^{1,0}M \oplus T^{0,1}M$, where $T_{\mathbb{C}}M$ is a complexification of $TM$, $$T^{1,0}M=\{X-iJX:X\in TM\}$$ and $$T^{0,1}M=\{X+iJX:X\in TM\}=\{\bar\zeta\in T_{\mathbb{C}}M:\zeta\in T^{1,0}M\}.$$

Let $\mathcal{A}^k$ be the set of $k$-forms, i.e. the set of sections of $\bigwedge^k(T_{\mathbb{C}}M)^\star$ and let
$\mathcal{A}^{p,q}$ be the set of $(p,q)$-forms, i.e. the set of sections of $\bigwedge^p(T^{1,0}M)^\star\otimes_{(\mathbb{C})}~\bigwedge^q(T^{0,1}M)^\star$. Then we have a direct sum decomposition $\mathcal{A}^k=\bigoplus_{p+q=k}\mathcal{A}^{p,q}$. We denote the projections $\mathcal{A}^k\rightarrow\mathcal{A}^{p,q}$ by $\Pi^{p,q}$.

Let $d:\mathcal{A}^k\rightarrow\mathcal{A}^{k+1}$ be (the $\mathbb{C}$-linear extension of) the exterior differential. Then $d=\da+\db-\ta-\tb$, where $\partial:=\Pi^{p+1,q}\circ d$, $\bar\partial:=\Pi^{p,q+1}\circ d$, $\ta:=-\Pi^{p+2,q-1}\circ d$, $\tb:=-\Pi^{p-1,q+2}\circ d$ on $\mathcal{A}^{p,q}$. Note that $\ta$ and $\tb$ are operators of order $0$. Let $\omega$ be a $(p,q)$-form. We have following formulas (see \cite{p}):
$$\partial\omega(\zeta_1,\ldots,\zeta_{p+1},\bar\eta_1,\ldots,\bar\eta_q)=\sum_{k=1}^{p+1}(-1)^{k+1}\zeta_k\omega(\zeta_1,\ldots,\widehat \zeta_k,\ldots,\bar\eta_q)$$
$$+\sum_{1\leq k<l\leq p+1}(-1)^{k+l}\omega([\zeta_k,\zeta_l],\zeta_1,\ldots,\widehat \zeta_k,\ldots,\widehat\zeta_l,\ldots,\bar\eta_q)$$ $$+\sum_{\begin{subarray}{c}
 1\leq k\leq p+1\\1\leq l\leq q
\end{subarray}
}(-1)^{k+l+p+1}\omega([\zeta_k,\bar\eta_l],\zeta_1,\ldots,\widehat \zeta_k,\ldots,\widehat{\bar\eta_l},\ldots,\bar\eta_q),$$
$$\bar\partial\omega(\zeta_1,\ldots,\zeta_{p},\bar\eta_1,\ldots,\bar\eta_{q+1})=\sum_{k=1}^{q+1}(-1)^{k+p+1}\bar\eta_k\omega(\zeta_1,\ldots,\widehat{\bar\eta_k},\ldots,\bar\eta_{q+1})$$
$$+\sum_{1\leq k<l\leq q+1}(-1)^{k+l}\omega([\bar\eta_k,\bar\eta_l],\zeta_1,\ldots,\widehat{ \bar\eta_k},\ldots,\widehat{\bar\eta_l},\ldots,\bar\eta_{q+1})$$ $$+\sum_{\begin{subarray}{c}
 1\leq k\leq p\\1\leq l\leq q+1
\end{subarray}
}(-1)^{k+l+p}\omega([\zeta_k,\bar\eta_l],\zeta_1,\ldots,\widehat \zeta_k,\ldots,\widehat{\bar\eta_l},\ldots,\bar\eta_q),$$
$$\ta\omega(\zeta_1,\ldots,\zeta_{p+2},\bar\eta_1,\ldots,\bar\eta_{q-1})$$ $$=-\sum_{1\leq k<l\leq p+2}(-1)^{k+l}\omega([\zeta_k,\zeta_l],\zeta_1,\ldots,\widehat \zeta_k,\ldots,\widehat\zeta_l,\ldots,\bar\eta_{q-1}),$$
$$\tb\omega(\zeta_1,\ldots,\zeta_{p-1},\bar\eta_1,\ldots,\bar\eta_{q+2})$$ $$=-\sum_{1\leq k<l\leq q+2}(-1)^{k+l}\omega([\bar\eta_k,\bar\eta_l],\zeta_1,\ldots,\widehat {\bar\eta_k},\ldots,\widehat{\bar\eta_l},\ldots,\bar\eta_{q+2}),$$
where $\zeta_1,\ldots,\zeta_{p+2},\eta_1,\ldots,\eta_{q+2}$ are vector fields of type $(1,0)$ (i.e. sections of $T^{1,0}M$).
  In particular, for a smooth function $u$ we have :
 $$i\partial\bar\partial u=i\sum (\zeta_p\bar{\zeta_q}u-[\zeta_p,\bar{\zeta_q}]^{0,1}u)\zeta_p^\star\wedge \bar\zeta_q^\star,$$
  where $\zeta_1,\ldots,\zeta_n$ is  a (local) frame of $T^{1,0}$ and $\zeta_1^\star,\ldots, \zeta_n^\star,\bar\zeta_1^\star,\ldots, \bar\zeta_n^\star$ is a base of $(T_{\mathbb{C}}M)^\star$ dual to the base $\zeta_1,\ldots,\zeta_n,\bar{\zeta_1},\ldots,\bar{\zeta_n}$ of $T_{\mathbb{C}}M$. We will use the following identities: $$\da\db+\db\da+\ta\tb+\tb\ta=0,\;\;\;\;\;\da^2=\ta\db+\db\ta,\;\;\;\;\db^2=\tb\da+\da\tb.$$
  
  We say that an almost complex structure $J$ is integrable, if any of the following (equivalent) conditions is satisfied:\\i) $d=\partial+\bar\partial$;
  \\ii) $\bar\partial^2=0$;
  \\iii) $[\zeta,\xi]\in T^{0,1}M$ for vector fields $\zeta,\xi\in T^{0,1}M$.
  \\ By the Newlander-Nirenberg Theorem $J$ is integrable if and only if it is induced by a complex structure.

We can define the positivity of $(p,p)$-forms or more general of $(p,p)$-currents, in the same way as on complex manifolds. Positive currents are of order $0$ (see \cite{p}).

We  always assume that there is a fixed hermitian metric $\omega$ on $M$, i.e.  $(1,1)$-positive form, and $L^p$ and $W^{1,2}$ norms and a distance on $M$ are defined with respect to this metric.

Let $\mathbb{D}=\{z\in\mathbb{C}:|z|<1\}$. We say that a (smooth) function $\lambda:\mathbb{D}\rightarrow M$ is $J$-holomorphic or simpler holomorphic, if $\lambda'(\frac{\partial}{\partial\bar z})\in T^{0,1}M$. The following proposition (see \cite{i-r1}) shows that there exist plenty of such disks:

\begin{prp}\label{dyski}
 Let $0\in M\subset\mathbb{R}^{2n}$, $k,k'\geq1$. For $v_0,v_1,\ldots,v_k\in\mathbb{R}^{2n}$ close enough to 0, there is a holomorphic function $\lambda:\mathbb{D}\rightarrow M$, such that $\lambda(0)=v_0$ and $\frac{\partial^l\lambda}{\partial x^l}=v_l$, for $l=1,\ldots,k$. Moreover, we can choose $\lambda$ with $\mathcal{C}^1$ dependence on parameters $(v_0,\ldots,v_k)\in(\mathbb{R}^{2n})^{k+1}$, where for holomorphic functions  we consider $\mathcal{C}^{k'}$ norm.
\end{prp}
If $\lambda:\mathbb{D}\rightarrow M$ is holomorphic and $u$ is a smooth function, then
 \begin{equation}\label{laplasjan}
\Delta (u\circ \lambda)=\frac{1}{2}\dc u(\zeta,\bar\zeta),\end{equation}
 where $\zeta=\frac{\partial\lambda}{\partial x}-iJ\frac{\partial\lambda}{\partial x}$ .

\subsection{plurisubharmonic functions}
An upper semi-continuous function $u$ on an open subset of $M$ is said to be $J$-plurisubharmonic or simpler plurisubharmonic, if a function $u\circ\lambda$ is subharmonic for every holomorphic function $\lambda$. We denote the set of plurisubharmonic functions on $\Omega\subset M$ by $\mathcal{PSH}(\Omega)$. Similarly as in the  case of integrable complex structure, plurisubharmonic functions are in spaces $L^p_{loc}$ for any $p<+\infty$. Recently Harvey and Lawson proved that an upper semicontinuous function localy integrable $u$ is plurisubharmonic iff a current $i\partial\bar\partial u$ is nonnegative (see \cite{h-l}). We say that a function $u $ on $\Omega$ is strictly plurisubharmonic if  for every open set $D\Subset\Omega$ and  $\mathcal{C}^2$ function $\varphi$ on a neighbourhood of $\bar D$ there is $\varepsilon>0$ such that a function $u+\varepsilon\varphi$ is plurisubharmonic in $D$.

We say that a domain $\Omega\Subset M$ is strictly pseudoconvex of class $\mathcal{C}^{\infty}$ (respectively of class $\mathcal{C}^{1,1}$),  if there is a strictly plurisubharmonic function $\rho$ of class $\mathcal{C}^{\infty}$ (respectively of class $\mathcal{C}^{1,1}$)  in a neighbourhood of $\bar\Omega$,  such that $\Omega=\{\rho<0\}$ and $\triangledown\rho\neq0$ on $\partial\Omega$. In that case we say that $\rho$ is a defining function for $\Omega$. 
\subsection{Monge-Amp\`ere equation}

The following is the main theorem in \cite{pl}:
\begin{thr}\label{MAtheorem}
Let $\Omega\Subset M$ be a strictly pseudoconvex domain of class $\mathcal{C}^{\infty}$. There is a unique  solution $u$ of the Dirichlet problem:
 \begin{equation*}\label{DP}
\left\{
\begin{array}{l}
    u\in\mathcal{PSH}(\Omega)\cap \mathcal{C}^\infty(\bar\Omega)\\ 
    (i\partial\bar\partial u)^n=dV \;\mbox{ in }\;\Omega\\
    u=\varphi\;\mbox{ on }\;\partial\Omega
\end{array}
\right.\;,
\end{equation*}
where $\varphi\in\mathcal{C}^\infty(\bar\Omega)$ and $dV$ is the volume form on a neighbourhood of $\bar\Omega$. 
\end{thr}

The following version of the comparison principle is also proved in the same paper:
\begin{prp}\label{cp'}
Suppose that   $u$, $v\in\mathcal{C}^2(\bar\Omega)$ are such that $v$ is a plurisubharmonic function and $(i\partial\bar{\partial} u)^n\leq(i\partial\bar{\partial} v)^n$ on the set $\{i\partial\bar{\partial} u>0\}$. Then for any $H\in\mathcal{PSH}$, an inequality $$\varlimsup_{z\rightarrow z_0}(v+H-u)\leq0$$  for any $z_0\in\po$, implies $v+H\leq u$ on $\Omega$.
\end{prp}

Harvey and Lawson solved the Dirichlet problem with continuous date:
\begin{thr}[see \cite{h-l}]\label{HLMAtheorem}
 Let $\Omega\Subset M$ be a strictly pseudoconvex domain of class $\mathcal{C}^{1,1}$. There is a unique viscosity solution $u$ of the Dirichlet problem:
 \begin{equation}\label{DP}
\left\{
\begin{array}{l}
    u\in\mathcal{PSH}(\Omega)\cap \mathcal{C}(\bar\Omega)\\ 
    (i\partial\bar\partial u)^n=fdV \;\mbox{ in }\;\Omega\\
    u=\varphi\;\mbox{ on }\;\partial\Omega
\end{array}
\right.\;,
\end{equation}
where $\varphi,f\in\mathcal{C}(\bar\Omega)$, $f\geq0$ and $dV$ is the volume form on a neighbourhood of $\bar\Omega$. 
\end{thr}

We can easily obtain the existence part of the above result from Theorem \ref{MAtheorem} (see \cite{pl} for details in case $f=0$, the general case can be proved almost in the same way). In particular the solution is a limit of smooth solutions of Dirichlet problems with a smooth date. Further, using the gradient estimate (Lemma 3.3 in \cite{pl}) we obtain the existence of the Lipschitz solution. 

\begin{thr}\label{LMAtheorem}
Let $\Omega\Subset M$ be a strictly pseudoconvex domain of class $\mathcal{C}^{1,1}$. If $\varphi\in\mathcal{C}^{1,1}(\bar\Omega)$ and $f^{1/n}\in\mathcal{C}^{0,1}(\bar\Omega)$ then there is a unique viscosity Lipschitz solution $u$ of the Dirichlet problem (\ref{DP}). 
\end{thr}

We will see in section \ref{Monge-Amp\`ere operator} that in case $n=2$ a viscosity solution of the Dirichlet problem is also a solution in the pluripotential sense.

\section{approximation}\label{approximation}

In this section we prove Richberg Theorem for plurisubharmonic functions on almost complex manifolds.

\begin{thr}\label{aproksymacjaciaglych}
If $u,h\in\mathcal{C}(M)$, $h>0$ and $u$ is strictly plurisubharmonic, then  there exists a strictly plurisubharmonic function $\psi\in\mathcal{C}^\infty(M)$ such that $u\leq\psi\leq u+h$.
\end{thr}
As an immediate consequence we obtain the following:

\begin{crl}
 If there is a continuous strictly plurisubharmonic function on $M$, there is also a smooth one.
\end{crl}

Theorem \ref{aproksymacjaciaglych} on complex manifolds was proved in \cite{r}. Non integrable case was stated in \cite{c-e} as an open problem. The corollary answers to the problem of Ivashkovich and Rosay from \cite{i-r2}, however they not assume $\mathcal{C}^\infty$ regularity of $J$, so the question ``Does  the existence of a continuous strictly J-plurisubharmonic function ensures the existence of the smooth one?'' remains open for non smooth $J$.

 We will need the following lemma:
\begin{lm}\label{lokalnerozwiazaniePD}
Let  $u$ be as in Theorem \ref{aproksymacjaciaglych}. If $ U\Subset M$ is a smooth strictly pseudoconvex domain 
 and  $K\Subset U$, then there is $v\in\mathcal{C}^\infty(\bar U)$  strictly plurisubharmonic on $U$ such that $v<u$ on $\partial U$ and $v>u$ on $K$.
\end{lm}

{\it Proof:} Let $\rho$ be a defining function for $U$ and let $\varepsilon>0$ be such that $u-\varepsilon\rho\in\mathcal{PSH}(U)$. We can take a  function $\varphi\in\mathcal{C}^\infty(\bar U)$ such that $$u+\frac{\varepsilon}{4}\sup_K\rho<\varphi<u.$$ Let $v$ be a solution of the following Dirichlet Problem 
$$
\left\{
\begin{array}{l}
    v\in\mathcal{PSH}(U)\cap \mathcal{C}^\infty(\bar U)\\ 
    (i\partial\bar\partial v)^n=(\frac{\varepsilon}{2})^n(i\partial\bar{\partial}\rho)^n \;\mbox{ in }U\\
    v=\varphi\;\mbox{ on }\;\partial U
\end{array}
\right.\;,
$$
then by the Comparison Principle $$v\geq u-\frac{\varepsilon\rho}{2}+\frac{\varepsilon}{4}\sup_K\rho>u \hbox{ on } K.\;\;\Box$$

For the proof of Theorem \ref{aproksymacjaciaglych} we need the regularised maximum of two functions. Let $m_s$ be a smooth convex function in $\mathbb{R}^2$, such that $\max\leq m_s\leq\max+s$ and $m_s(x,y)=\max\{x,y\}$, if $|x-y|\geq s$. If $U_1,U_2\subset M$ and $u_i$ is a function on $U_i$, we can put $${\max}_s\{u_1,u_2\}=\left\{
\begin{array}{l}
    m_s(u_1,u_2) \mbox{ on } U_1\cap U_2\\ 
    u_1 \mbox{ on } U_1\setminus U_2\\
u_2 \mbox{ on } U_2\setminus U_1
\end{array}
\right.\;. $$ Obviously if $u_1$, $u_2$ are plurisubharmonic, $u_1+s<u_2$ on $\partial U_1\cap U_2$ and $u_2+s<u_1$ on $\partial U_2\cap U_1$, then the function ${\max}_s\{u_1,u_2\}$ is plurisubharmonic on $U_1\cup  U_2$. The following Lemma is also  proved by Sukhov (see \cite{s}):

\begin{lm}\label{max}
 Let $u_1$, $u_2$ be as above. If additionally they are smooth strictly plurisubharmonic, then $\max_s\{u_1,u_2\}$ is also (smooth) strictly plurisubharmonic.
\end{lm}

{\it Proof:} If $\max_s\{u_1,u_2\}$ is not strictly plurisubharmonic, then there are a point $z_0\in U_1\cup  U_2$ and a vector $\zeta=X-iJX\in T_{z_0}^{1,0}$, such that $\dc\max_s\{u_1,u_2\}(z_0)(\zeta,\bar\zeta)=0$. By Proposition \ref{dyski} there is a holomorphic disc $\lambda:\mathbb{D}\rightarrow M$ such that $\lambda(0)=z_0$ and $\frac{\da}{\da x}(0)=X$. Functions $u_1\circ\lambda$, $u_2\circ\lambda$ are strictly plurisubharmonic and it is easy to calculate that $\Delta\max_s\{u_1\circ\lambda, u_2\circ\lambda\}(0)>0$. Using formula (\ref{laplasjan}) we get the contradiction. $\;\Box$

{\it Proof of Theorem \ref{aproksymacjaciaglych}:} For every $z\in M$ we can choose open sets $U,V$ such that $z\in V\Subset U\Subset M$, $U$ is strictly pseudoconvex and $\sup_U u<\inf_U (u+h)$. Hence there  are $(V_n)_{n=1}^\infty$, $(U_n)_{n=1}^\infty$ locally finite open covers of $M$, which for every $n$ satisfy the following conditions:
 \\1) $V_n\Subset U_n$;
 \\2) $U_n$ is strictly pseudoconvex;
 \\3) $\sup_{U_n} u<\inf_{U_n} (u+h)$. \\
Let $W_0=\emptyset$ and $W_n=\bigcup_{k=1}^nV_n$. By Lemma \ref{lokalnerozwiazaniePD} for every $n$ there is a  strictly plurisubharmonic function   $v_n\in\mathcal{C}^\infty(\overline{U_n})$, such that $v_n<u$ on $\partial U_n$ and $v_n>u$ on $\overline{V_n}$. Note that $v_n<u+h$.

Let $$S_p=\frac{\inf_{\overline{ V_{p}}}(v_p-u)}{\#\{k:U_k\cap U_p\neq\emptyset\}}$$ and let $$S^n_p=\#\{k\leq n:U_k\cap U_p\neq\emptyset\}S_p.$$ Let us define, by the induction, a sequence  $\psi_n$ of continuous plurisubharmonic functions on $M$, which  satisfy the following conditions: \\
i) $u\leq\psi_n<u+h$,\\
ii) $\psi_n>u$ on $\overline{W_n}$,\\
iii) $\psi_n$ is strictly plurisubharmonic and smooth in $W_n$,\\
iv) $\psi_n$ is strictly plurisubharmonic and smooth also in all sets\\ $\{~\psi_n~>~u~+~S_p^n~\}\cap~U_p$.\\
 Let $\psi_0=v_0$. Now assume that $\psi_n$ is as above. We can choose $s>0$, such that  $$s<\inf_{U_{n+1}}(u+h-\max\{v_{n+1},\psi_n\}),$$ $$s<\inf_{\partial U_{n+1}}(\psi_n-v_{n+1}),$$
and for every $p$, such that $U_{n+1}\cap U_p\neq\emptyset$, we have $s<S_p$.
 Then we can put $\psi_{n+1}=\max_s\{v_{n+1},\psi_n\}$. Obviously $\psi_{n+1}$ satisfies above conditions i) and ii). Note that $\psi_{n+1}$ is strictly plurisubharmonic and smooth on the sets $\{\psi_{n+1}>\psi_n+S_p\}$, hence $\psi_{n+1}$ satisfies iv) and so iii).
 
 Observe that for any compact set $K\subset M$, there is $n_0\in\mathbb{N}$ such that $\psi_n=\psi_{n_0}$ on $K$ for $n\geq n_0$. Therefore $$\psi=\lim_{n\rightarrow\infty}\psi_n$$ is as in the Theorem. $\;\Box$

The idea of using solutions of the Monge-Amp\`{e}re equation to approximate plurisubharmonic functions is probably due to J. P. Rosay.

\section{$W^{1,2}$ estimates for plurisubharmonic functions}

In this section we prove some properties of plurisubharmonic functions in the Sobolew space $W^{1,2}_{loc}$.

The following lemma is the special case of Theorem 3.3  from \cite{b1}.
\begin{lm}\label{blockiw12} Let $D\Subset\mathbb{D}$, $u,v\in\mathcal{SH}$, $u\leq v\leq0$ and $u\in{W}^{1,2}(\mathbb{D})$. Then
 $v\in{W}^{1,2}(D)$ and $\|v\|_{W^{1,2}(D)}\leq C\|u\|_{W^{1,2}(\mathbb{D})}$, where the constant $C$ depends only on $D$.
\end{lm}

\begin{prp}\label{W12estimates} Let $u\in\mathcal{PSH}\cap{W}^{1,2}_{loc}(\Omega)$ then:\\
i) If $v\in\mathcal{PSH}(\Omega)$ and $v\geq u$,  then $v\in{W}^{1,2}_{loc}(\Omega)$;\\
ii) If a sequence $u_j$ of plurisubharmonic functions decreases to $u$, then it converges in ${W}^{1,2}_{loc}$;\\
iii) If a sequence $u_j$ of bounded plurisubharmonic functions increases to $u$ a. e., then it converges in ${W}^{1,2}_{loc}$.
\end{prp}

Corresponding results for subharmonic functions in $\mathbb{R}^m$ were proved in \cite{b1} (first part) and \cite{c} (second part). Błocki gives a nice proof of both in \cite{b2} and his proof of $ii)$ works also in our case for $ii)$ and $iii)$ .

{\it Proof:} We can assume that $u\in\mathcal{PSH}\cap{W}^{1,2}(\Omega)$. To prove $i)$ we can assume $v\leq0$. Let $z_0\in\Omega$.
By Proposition \ref{dyski} we can choose $\mathcal{C}^\infty$ embeddings  $$H_k:B\times\bar{\mathbb{D}}\rightarrow \Omega,
 \;\;k=1,\ldots,n,$$
where $B=\{z\in\mathbb{C}^{n-1}:|z|\leq1\}$, such that \\1) $H_k(0,0)=z_0$,\\2)  functions $$h_{k,t}:\bar{\mathbb{D}}\rightarrow M, \;\;k=1,\ldots,n,\;t\in B,$$ given by $h_{k,t}(z)=H_k(t,z)$, are holomorphic,\\3)  for every $w\in D=\bigcap_kH_k(\int B\times \mathbb{D})$ and points $(t_k,z_k)=H_k^{-1}(w)$, vectors 
$\frac{\partial}{\partial x}h_{k,t_k}(z_k)$ are $J$ linear independent in $T_{w}M$.

Let us put $u_{k,t}=u\circ h_{k,t}$, $v_{k,t}=v\circ h_{k,t}$ and let $$U=\{z\in\mathbb{C}:|z|<\frac{1}{2}\}\Subset\mathbb{D}.$$ By 
Lemma \ref{blockiw12} we obtain that
 $$\int_B\left(\int_U i\da v_{k,t}\wedge\db v_{k,t}\right)d\lambda\leq C\int_B\left(\int_\mathbb{D}\left( i\da u_{k,t}\wedge\db u_{k,t}+|u|^2i\da z\wedge\da\bar z\right)\right)d\lambda,$$ where $d\lambda$ is $(2n-2)$-dimensional Lebesgue measure. To conclude that 
\begin{equation}\label{w12dlamniejszej}
 \int_{D'}i\da v\wedge\db v\wedge\omega^{n-1}\leq C(J,\omega)\|u\|_{W^{1,2}(\Omega)}^2<+\infty,
\end{equation}
where $D'=\bigcap_kH_k(\int B\times U)$ is a neighbourhood of $z_0$, it is enough to note that there is a constant $C_1$ (depending on functions $H_k$), such that
$$\int_B\left(\int_\mathbb{D}\left( i\da u_{k,t}\wedge\db u_{k,t}+|u|^2i\da z\wedge\da\bar z\right)\right)d\lambda\leq C_1\|u\|_{W^{1,2}(\Omega)}^2$$
for $k=1,\ldots,n$ and
$$\int_{D'}i\da v\wedge\db v\wedge\omega^{n-1}\leq C_1\sum_k\int_B\left(\int_U i\da v_{k,t}\wedge\db v_{k,t}\right)d\lambda.$$

To prove $ii)$ and $iii)$ note that by $i)$ all $u_j$ are in ${W}^{1,2}_{loc}(\Omega)$. We need the following fact: $v_1\dc v_2$ is well defined current of order $0$ for $v_1,v_2\in\mathcal{PSH}\cap{W}^{1,2}_{loc}(\Omega)$. Let $\varphi$ be a non negative smooth function with a compact support in $\Omega$ and $V_k=\max\{v_1,-k\}$. Note that by the estimation (\ref{w12dlamniejszej}) a $W^{1,2}$ norm of $V_k$, on a set $\{\varphi>0\}$, does not depend on $k$. Using Stokes' theorem we can estimate $$\left|\int_\Omega\varphi V_k\dc v_2\wedge\omega^{n-1}\right|=\left|\int_\Omega\left(V_ki\db v_2\wedge\da(\varphi\omega^{n-1})-\varphi i\da V_k\wedge\db v_2\wedge\omega^{n-1}\right)\right| $$ $$\leq C(\varphi)\|V_k\|_{W^{1,2}(\{\varphi>0\})}\|v_2\|_{W^{1,2}(\{\varphi>0\})}\leq C(\varphi, v_1,v_2)$$
and thus we get  $v_1\in L_{loc}^1(\Omega,\dc v_2\wedge\omega^{n-1})$.

 Now, as in \cite{b2}, we can take a non negative smooth function $\varphi$ with a compact support in $\Omega$ and (again) using Stokes' theorem we get: 
$$\int_\Omega\varphi i\da(u_j-u)\wedge\db(u_j-u)\wedge\omega^{n-1}$$
$$=\int_\Omega\varphi\left( \frac{1}{2}\dc(u_j-u)^2-(u_j-u)\dc(u_j-u)\right)\wedge\omega^{n-1}$$
$$=\int_\Omega \frac{1}{2}(u_j-u)^2\dc(\varphi\omega^{n-1})-\int_\Omega(u_j-u)\dc(u_j-u)\wedge\omega^{n-1}.$$
The first integral tends to $0$ by the Lebesgue monotone convergence theorem. The current $\dc (u_j-u)$ converges weakly to $0$ and this gives us that the second integral converges to 0 if $u_j$ and $u$ are bounded. If a sequence $u_j$ is decreasing  then we can estimate the second integral:
$$-\int_\Omega(u_j-u)\dc(u_j-u)\wedge\omega^{n-1}\leq \int_\Omega(u_j-u)\dc u\wedge\omega^{n-1}$$ and the last integral (again by the Lebesgue monotone convergence theorem) converges to 0. $\;\Box$

Let us see again on the proof of $i)$. If $u$ is a constant function and we decrease proportionally $\Omega$, $D$ and $D'$, then the constant $C(J,\omega)$ in (\ref{w12dlamniejszej}) is decreasing in a controlled way too. In particullar we obtain the following:

\begin{prp}\label{W12estimates'} Let $z_0\in M'\Subset M$. Let $\Omega(R)=\{z:{\rm dist}(z,z_0)<R\}$ and $D(R)=\{z:{\rm dist}(z,z_0)<\frac{R}{2}\}$. Then there are $R_0>0$ and $C=C(M',\omega,J)$, such that for every $R<R_0$ and $v\in\mathcal{PSH}\cap L^\infty(\Omega(R))$ we have $$\|v\|_{W^{1,2}(D(R))}\leq C R^{2n-2}\|v\|_{L^\infty(\Omega(R))} .$$
\end{prp}

\section{Monge-Amp\`ere operator}\label{Monge-Amp\`ere operator}
\subsection{definition}\label{definition}
Let $u,v\in\mathcal{PSH}\cap W^{1,2}_{loc}(\Omega)$ then 
 \begin{equation*}\label{wedgeproduct}
\dc u\wedge\dc v:=-\dc(i\da u\wedge\db v)+\da(\da u\wedge\tb\da v)+\db(\ta\db u\wedge\db v)+\tc\da u\wedge\db v-\ta\db u\wedge\tb\da v
\end{equation*}
is a well defined $(2,2)$ current. If $u$ or $v$ can be approximated in $W^{1,2}_{loc}(\Omega)$ by smooth plurisubharmonic functions, then this is a positive current,  moreover if $n=2$ this is a positive measure. Note that if $u$ and $v$ are of class $\mathcal{C}^{1,1}$, it is a usual wedge product of (bounded) forms.

Let $$\mathcal{D}:=\{u\in\mathcal{PSH}\cap W^{1,2}_{loc}(\Omega):\hbox{ a  current }(\dc u)^2 \hbox{ is positive}\}.$$
The following theorem is an easy consequence of Theorem \ref{aproksymacjaciaglych} and Proposition \ref{W12estimates}.\begin{thr}\label{MAdlaciaglych}
$\mathcal{PSH}\cap \mathcal{C}(\Omega)\subset \mathcal{D}$.
\end{thr}

Proposition \ref{W12estimates} gives us the following
\begin{prp}\label{zbieznoscMA}
Let $u\in\mathcal{PSH}\cap{W}^{1,2}_{loc}(\Omega)$ then:\\
i) If a sequence $(u_j)\subset \mathcal{D}$  decreases to $u$, then $u\in\mathcal{D}$ and $(\dc u_j)^2$ converges to $(\dc u)^2$;\\
ii) If a sequence $(u_j)\subset \mathcal{D}\cap L^\infty_{loc}$  increases to $u$ a.e., then $u\in\mathcal{D}$ and $(\dc u_j)^2$ converges to $(\dc u)^2$.
\end{prp}

In $\mathbb{C}^2$ (with the standard almost complex structure) the  operator $\ma$, was defined by Bedford and Taylor , for $W^{1,2}_{loc}$ plurisubharmonic functions   (essensialy as $-\dc(i\da u\wedge\db v)$, see \cite{b-t2}). In \cite{b1} Błocki proved  that $W^{1,2}_{loc}\cap\mathcal{PSH}$  forms a  natural domain for the Monge-Amp\`ere operator.

The following example shows that there are also unbounded functions for which Monge-Amp\`ere operator is well defined.$\newline\newline$
{\it Example:} For every point in $M$ there is a plurisubharmonic function (on some neighbourhood of this point) with a logarithmic singularity in this point. Let $J$ be an almost complex structure in $\mathbb{C}^2$ such that $J(0)=J_{st}(0)$. Then there is $A>0$, such that a function $L(z)=\log|z|+A|z|$ is plurisubharmonic in some neighbourhood $U$ of $0$. Such functions are crucially used  in order to localize and estimate the Kobayashi-Royden metric on an almost complex manifold (see for example \cite{g-s} and \cite{b}).

A function $L$ is in $W^{1,2}_{loc}(U)$ and a sequence of continuous functions $\max\{-k,L\}$ decreases to $L$. Therefore the Monge-Amp\`ere measure $(\dc L)^2$ is well defined. Outside $0$ it is a smooth volume form. Now we calculate $(\dc L)^2(\{0\})$. Note that $(\dc |z|^2)^2=fdV$, where $f$ is a smooth function, $f(0)=8$ and $dV$ is the standard volume form in $\mathbb{C}^2$. Let us put $$L_k=\left\{ \begin{array}{lll}
\frac{1}{2}k(k+A)|z|^2-\frac{k+3A}{2k}-\log k & \mbox{ on } & B_k=\{|z|\leq \frac{1}{k}\}\\
L & \mbox{ on }  & U\setminus B_k.
\end{array}\right.
$$
For $k$ large enough $L_k$ is $\mathcal{C}^{1,1}$ plurisubharmonic function on $U$, $L_k\searrow L$ and $$\frac{(k+A)^2}{8k^2}\pi^2\min_{B_k}f\leq(\dc L_k)^2(B_k)\leq \frac{(k+A)^2}{8k^2}\pi^2\max_{B_k}f.$$ We can conclude that $(\dc L)^2(\{0\})=\pi^2$.

\subsection{current $(\dc)^2u$}
From now we will always assume that $n=2$.

In this subsection we consider an operator $(\dc)^2$, which will appear naturally in the proof of the comparison principle. Let us calculate:
$$ (\dc)^2 =\ta\db\tb\da+\db\tc\da+\da\tb\ta\db\stackrel{n=2}{=}\db\tc\da+\da\tb\ta\db=\db\tc\da+\overline{\db\tc\da}.$$
In particular $(\dc)^2$ is a real operator. Let $u\in \mathcal{C}^2(\Omega)$. 
$$\db\tc\da u(\zeta_1,\zeta_2,\bar\zeta_1,\bar\zeta_2)=\bar\zeta_1(\tc\da u(\zeta_1,\zeta_2,\bar\zeta_2))-\bar\zeta_2(\tc\da u(\zeta_1,\zeta_2,\bar\zeta_1))$$
$$-\tc\da u(\zeta_1,\zeta_2,[\bar\zeta_1,\bar\zeta_2]^{0,1})$$ $$+\tc\da u([\zeta_1,\bar\zeta_1]^{1,0},\zeta_2,\bar\zeta_2)-\tc\da u(\zeta_1,[\zeta_2,\bar\zeta_2]^{1,0},\bar\zeta_1)$$ $$- \tc\da u([\zeta_1,\bar\zeta_2]^{1,0},\zeta_2,\bar\zeta_1) +\tc\da u(\zeta_1,[\zeta_2,\bar\zeta_1]^{1,0},\bar\zeta_2)$$

$$\bar\zeta_1(\tc\da u(\zeta_1,\zeta_2,\bar\zeta_2))=\bar\zeta_1\left[[\zeta_1,\zeta_2]^{0,1},\bar\zeta_2\right]^{1,0}u$$ $$=\da\db u(\left[[\zeta_1,\zeta_2]^{0,1},\bar\zeta_2\right]^{1,0},\bar\zeta_1)-\left[\left[[\zeta_1,\zeta_2]^{0,1},\bar\zeta_2\right]^{1,0},\bar\zeta_1\right]^{1,0}u,$$

$$\bar\zeta_2(\tc\da u(\zeta_1,\zeta_2,\bar\zeta_1))=\bar\zeta_2\left[[\zeta_1,\zeta_2]^{0,1},\bar\zeta_1\right]^{1,0}u$$ $$=\da\db u(\left[[\zeta_1,\zeta_2]^{0,1},\bar\zeta_1\right]^{1,0},\bar\zeta_2)-\left[\left[[\zeta_1,\zeta_2]^{0,1},\bar\zeta_1\right]^{1,0},\bar\zeta_2\right]^{1,0}u$$
and we can conclude that
$$(\dc)^2u(\zeta_1,\zeta_2,\bar\zeta_1,\bar\zeta_2)$$ $$=T_Ju+2{\rm Re}\left(\da\db u(\left[[\zeta_1,\zeta_2]^{0,1},\bar\zeta_2\right]^{1,0},\bar\zeta_1)-\da\db u(\left[[\zeta_1,\zeta_2]^{0,1},\bar\zeta_1\right]^{1,0},\bar\zeta_2)\right),$$
where $T_J$ is a real vector field.  The above formula gives us the following lemma:

\begin{lm}\label{oszacowaniedckwadrat}
Let $u$ be as above. Then there is a constant $C$, which depends only on $\Omega$, $J$ and $\omega$, such that 
$$-C\left(|T_Ju|\omega^2+\dc u\wedge\omega\right)\leq(\dc)^2u\leq C\left(|T_Ju|\omega^2+\dc u\wedge\omega\right).$$
In particular $(\dc)^2u$ is of order $0$.
\end{lm}

Note that $(\dc)^2u(\zeta_1,\zeta_2,\bar\zeta_1,\bar\zeta_2)$ and $$\left(\da\db u(\left[[\zeta_1,\zeta_2]^{0,1},\bar\zeta_2\right]^{1,0},\bar\zeta_1)-\da\db u(\left[[\zeta_1,\zeta_2]^{0,1},\bar\zeta_1\right]^{1,0},\bar\zeta_2)\right)$$ depends only on $J$ and on the volume form $\zeta_1^\star\wedge\zeta_2^\star\wedge\bar\zeta_1^\star\wedge\bar\zeta_2^\star$, so $T_J$ depends only on them as well. It seems to be interesting that on a manifold $(N=\{z\in M:T_J(z)\neq 0\},J|_N)$ we get a canonical (real) line bundle $L=\{xT_J:x\in\mathbb{R}\}$ with a canonical orientation.

The next example shows that it is possible that $T_J= 0$, even if $J$ is not integrable (and it suggests that it is not a generic case).$\newline\newline$
{\it Example:}  Let $(x_1,y_1,x_2,y_2)$ be (real) coordinates of $\mathbb{C}^2$, and $a:\mathbb{C}^2\rightarrow \mathbb{R}$ be a smooth function.
Let $J_a$ be an almost complex structure with the following matrix representation:
$$J_a=\left[\begin{array}{cccc}
    0 & 1 & 0 & 0 \\ 
    -1 & 0 & 0 & 0 \\
    0 & 0 & a & 1 \\
    0 & 0 & -1-a^2 & -a 
\end{array}
\right]\;.$$
Let $\zeta_1=\frac{\da}{\da x_1}-iJ_a\frac{\da}{\da x_1}$ and $\zeta_2=\frac{\da}{\da x_2}-iJ_a\frac{\da}{\da x_2}$. Then $[\zeta_1,\bar\zeta_1]=0$,
$$[\zeta_1,\zeta_2]=-[\zeta_1,\bar\zeta_2]=\alpha\zeta_2+\beta\bar\zeta_2,$$
$$[\zeta_2,\bar\zeta_2]=\gamma\zeta_2+\delta\bar\zeta_2,$$
where
$$\alpha=\left(\frac{a}{1-ai}-\frac{i}{2}\right)\zeta_1a,\;\;\beta=-\left(\frac{a}{1+ai}+\frac{i}{2}\right)\zeta_1a,$$
$$\gamma=\left(\frac{a}{1-ai}-\frac{i}{2}\right)(\bar\zeta_2-\zeta_2)a,\;\;\delta=-\left(\frac{a}{1+ai}+\frac{i}{2}\right)(\bar\zeta_2-\zeta_2)a.$$
By the Newlander-Nirenberg Theorem $J_a$ is integrable if and only if $[\zeta_1,\zeta_2]^{0,1}=0$ on $\mathbb{C}^2$.  And we can see that it is exactly when $\zeta_1a=0$. 

For $J_a$ and $u\in \mathcal{C}^2(\Omega)$ one can compute:
$$\db\tc\da u(\zeta_1,\zeta_2,\bar\zeta_1,\bar\zeta_2)=-\bar\zeta_2(\tc\da u(\zeta_1,\zeta_2,\bar\zeta_1))-\tc\da u(\zeta_1,[\zeta_2,\bar\zeta_2]^{1,0},\bar\zeta_1)$$
$$=-i|\beta|^2\dc u(\zeta_2,\bar{\zeta_2})+(2\gamma|\beta|^2-\bar{\zeta_2}(|\beta|^2))\zeta_2u,$$ therefore $$(\dc)^2u(\zeta_1,\zeta_2,\bar\zeta_1,\bar\zeta_2)=T_{J_a}u=\left((2\gamma|\beta|^2-\bar{\zeta_2}(|\beta|^2))\zeta_2+(2 \bar\gamma|\beta|^2-{\zeta_2}(|\beta|^2))\bar\zeta_2\right)u,$$ 
thus $T_J$ is equal to $0$ iff $(\dc)^2=0$. If $a$ depends only on $x_1,y_2$, then $T_J=0$. On the other hand if $a$ depends also on $x_2,y_2$, it seems that  $T_J$ vanishes very rarely.

\subsection{comparison principle}

In the pluripotential theory on complex manifolds the comparison principle is a very effective tool. In particular it gives us the uniqueness of the solution of the Dirichlet Problem. In this subsection we prove some versions of the comparison principle in the non integrable case but with additional assumptions. In all propositions below $\Omega\subset M$ is a domain which attains a bounded continuous strictly plurisubharmonic function (by Proposition \ref{aproksymacjaciaglych} a domain $\Omega$ attains also a bounded smooth strictly plurisubharmonic function).  In the proofs below $C$ is a constant under control, but it can change from a line to a next line. 

 The following proposition shows that the comparison principle holds for Lipschitz plurisubharmonic functions and if $T_J=0$ for all continuous plurisubharmonic functions.

\begin{prp}
 Let $u$, $v\in\mathcal{PSH}(\Omega)$ be (locally) Lipschitz. If $(i\partial\bar\partial u)^2\leq(i\partial\bar\partial v)^2$ and $\varliminf_{z\rightarrow\partial\Omega}(u-v)\geq0$, then $v\leq u$.
\end{prp}

{\it Proof:} Assume that  a set $\{v>u\}$ is not empty. Thus, we can choose $\varepsilon>0$ and a negative smooth strictly plurisubharmonic function $\rho$, such that a set $\{v+\varepsilon\rho>u\}$ is also not empty. Let $z_0$ be a point where a function $v+\varepsilon\rho-u$ attains its maximum. Let us put $$w=v+\varepsilon\rho-\varepsilon'{\rm dist}(\cdot,z_0)^2+a,$$ where  $a,\varepsilon'>0$ are such that $w(z_0)=u(z_0)$ and a function $\chi=\varepsilon\rho-\varepsilon'{\rm dist}(\cdot,z_0)^2$ is strictly plurisubharmonic in some neighbourhood of $z_0$. Consider the set $$\Omega_s=\{u-w<s\}.$$ For small enough $s>0$ we have $\Omega_s\subset
\{{\rm dist}(\cdot,z_0)<\sqrt{\frac{s}{\varepsilon'}}
\}$ and $\chi$ is strictly plurisubharmonic in some neighbourhood of $\overline{\Omega_s}$. We can choose $s$ (as small as we want) such that $$(i\partial\bar\partial u)^2(\partial\Omega_s)=i\partial\bar\partial u\wedge\omega(\partial\Omega_s)=0.$$ Let $w_s=\max\{w+s,u\}$. Constants $C$ below do not depend on $s$. The integration by parts gives  us an estimate \begin{equation}\label{CP1}
\int_\Omega (\dc u\wedge\omega-\dc w_s\wedge\omega)=\int_\Omega (w_s-u)\dc\omega\leq Cs\lambda(\Omega_s).\end{equation}
We can easily estimate
\begin{equation}\label{MAjestduzy}\begin{array}{l}
 \int_\Omega \left((\dc w_s)^2-(\dc u)^2\right)\geq\frac{1}{2}\int_{\Omega_s} \left((\dc \chi)^2+\dc \chi\wedge\dc w_s\right)\\
\\
 \geq C^{-1}\left(\lambda(\Omega_s) +\int_{\Omega_s}\dc w_s\wedge\omega\right).
\end{array}
\end{equation}
On the other hand by the integration by parts and by (\ref{CP1}) we have
$$\int_\Omega \left((\dc w_s)^2-(\dc u)^2\right)=\int_\Omega ( u-w_s)\dc\dc(u+w_s) $$
 $$\leq C \int_\Omega\left((w_s-u)\dc(u+w_s)\wedge\omega+|(w_s-u)T_J(u+w_s)|\omega^2\right)$$
$$\leq Cs\left(\lambda(\Omega_s) +\int_{\Omega_s}\dc w_s\wedge\omega\right).$$
Choosing $s$ enough small we obtain the contradiction with (\ref{MAjestduzy}). $\;\Box$

In the similar way we can prove that the comparison principle holds among other for H\"older continuous plurisubharmonic functions.

\begin{prp}
 Let $u$, $v\in\mathcal{PSH}\cap\mathcal{C}(\Omega)$ satisfy $$|(u-v)(z)-(u-v)(z')|\leq\frac{c}{(\log({\rm dist}(z,z')))^4} $$ for some $c\in\mathbb{R}$ and $z,z'\in\Omega$ close enough. If $(i\partial\bar\partial u)^2\leq(i\partial\bar\partial v)^2$ and $\varliminf_{z\rightarrow\partial\Omega}(u-v)\geq0$, then $v\leq u$.
\end{prp}

{\it Proof:} We can use notation as above and assume that $M\subset\mathbb{C}^2$, $z_0=0$ and $J(0)=J_{st}(0)$. Let $\tilde w=w+s(2\log k)^{-1}L_k$, where $k~=~e^{\sqrt[4]{C'/s}}$, $L_k$ are functions  from Example in section 5.1 and $C'$ is chosen such that $$|(u-w)(z)-(u-w)(z')|\leq\frac{C'}{(\log|z-z'|)^4} $$ for $z,z'\in\Omega$ close enough. Consider the set $$\tilde\Omega_s=\{u-\tilde w<s\}.$$ For small enough $s>0$ we have $$B_k=\{|z|\leq \frac{1}{k}\}\subset\tilde\Omega_s\subset
\{|z|<\sqrt{\frac{2 s}{\varepsilon'}}
\},$$  functions $\chi$ and $L_k$ are  plurisubharmonic in the set  $\{|z|<\sqrt{2\frac{2 s}{\varepsilon'}}
\}$, and $(\dc L_k)^2(B_k)\geq 1$. We can choose $s$ such that $(i\partial\bar\partial u)^2(\partial\tilde\Omega_s)=i\partial\bar\partial u\wedge\omega(\partial\tilde\Omega_s)=0$. Let $\tilde w_s=\max\{\tilde w+s,u\}$. Similarly as before we get
\begin{equation*}
\int_\Omega (\dc u\wedge\omega-\dc\tilde w_s\wedge\omega)\leq Cs\lambda(\tilde\Omega_s),\end{equation*}
\begin{equation*}
 \int_\Omega \left((\dc\tilde w_s)^2-(\dc u)^2\right)
 \geq \left(s(2\log k)^{-1}\right)^2+C^{-1}\int_{\tilde\Omega_s}\dc\tilde w_s\wedge\omega
\end{equation*}
\begin{equation*}
\geq C^{-1}\left(s^2\sqrt s+\int_{\tilde\Omega_s}\dc\tilde w_s\wedge\omega\right)
\end{equation*}
and 
$$\int_\Omega \left((\dc\tilde w_s)^2-(\dc u)^2\right)\leq Cs \left(\int_{\tilde\Omega_s}\dc\tilde w_s\wedge\omega+\int_{\tilde\Omega_s}|T_J(u+\tilde w_s)|\omega^2\right).$$
To obtain the same contradiction, as in the previous proof, it is enough to estimate the last integral. Using Proposition \ref{W12estimates'} we get
$$ \int_{\tilde\Omega_s}|T_J(u+\tilde w_s)|\omega^2\leq C\|1\|_{W^{1,2}(\Omega_s)}\|u+\tilde w_s\|_{W^{1,2}(\Omega_s)}$$ $$\leq Cs\sqrt s\left(\log\left(\max_{|z|\leq2\sqrt{\frac{2 s}{\varepsilon'}}} (u(z)+\tilde w_s(z))-\min_{|z|\leq2\sqrt{\frac{2 s}{\varepsilon'}}} (u(z)+\tilde w_s(z))\right)\right)^{-1}.$$
Note that by the continuity, the expression inside the logarithm tends to 0, as $s$ tends to $0$. $\;\Box$

From the above proposition we get the following.
\begin{crl}
 There is at most one H\"older continuous solution of the following Dirichlet problem
\begin{equation*}\label{DP}
\left\{
\begin{array}{l}
    u\in\mathcal{PSH}(\Omega)\cap \mathcal{C}(\bar\Omega)\\ 
    (i\partial\bar\partial u)^2=\mu \;\mbox{ in }\;\Omega\\
    u=\varphi\;\mbox{ on }\;\partial\Omega
\end{array}
\right.\;,
\end{equation*}
where  $\mu$ is a Borel measure on $\Omega$ and $\varphi\in\mathcal{C}(\bar\Omega)$.
\end{crl}

It is not clear to the author, even for the Dirichlet problem with smooth date, how to prove the uniqueness of the solution in the class $\mathcal{D}$ or in $\mathcal{PSH}\cap\mathcal{C}(\Omega)$ .

$\newline$
\textbf{Acknowledgements.} The author would like to express his
gratitude to Z. Błocki for helpful discussions and advice during
the work on this paper.


\begin{thebibliography}{9999}
\bibitem[B-T1]{b-t1}
 E. Bedford, B. A. Taylor,
\textit{The Dirichlet problem for a complex Monge-Amp\`{e}re
 equation}, Invent. Math. 37(1976), 1-44,
\bibitem[B-T2]{b-t2} E. Bedford, B. A. Taylor,
\textit{Variational properties of the complex Monge-Amp\`{e}re equation
   I. Dirichlet principle}, Duke. Math. J. 45, 375-403 (1978),
 \bibitem[B-T3]{b-t3}
 E. Bedford, B. A. Taylor,
\textit{ A new capacity for plurisubharmonic functions}, Acta Math. 149(1982), 1-41,
\bibitem[B]{b}
 F. Bertrand, \textit{Sharp estimates of the Kobayashi metric and Gromov hyperbolicity},
J. Math. Anal. Appl. 345 (2008), no. 2, 825-844, 
\bibitem[B1]{b1}
Z. Błocki, \textit{On the definition of the Monge-Amp\`{e}re operator in $\mathbb{C}^2$}, Math.
Ann. 328 (2004),
415-423,
\bibitem[B2]{b2}
Z. Błocki, \textit{Minicourse on pluripotential theory}, University of Vienna, September 2012,

\bibitem[C]{c}
U. Cegrell, \textit{The gradient lemma}, Ann. Polon. Math. 91 (2007), 143-146,
\bibitem[C-E]{c-e}
K. Cieliebak, Y. Eliashberg, \textit{Stein Structures: existence and flexibility},
 ArXiv:1305.1619,
\bibitem[G-S]{g-s}
 H. Gaussier, A. Sukhov, \textit{Estimates of the Kobayashi-Royden metric in almost complex manifolds}, Bull. Soc. Math. France 133 (2005), no. 2, 259-273,
\bibitem[H-L]{h-l}
R. Harvey, B. Lawson,
\textit{Potential Theory on Almost Complex Manifolds}, ArXiv:1107.2584,
\bibitem[I-R1]{i-r1}
   S. Ivashkovich,  J.-P. Rosay,
   \textit{Schwarz-type lemmas for solutions of $\overline\partial$-inequalities and complete hyperbolicity of almost complex manifolds},
Ann. Inst. Fourier (Grenoble) 54 (2004), no. 7, 2387-2435,
 \bibitem[I-R2]{i-r2}
   S. Ivashkovich,  J.-P. Rosay,  \textit{Boundary values and boundary uniqueness of J-holomorphic mappings}, Int. Math. Res. Not. IMRN 2011, no. 17, 3839-3857,
\bibitem[K]{k} U. Kuzman,  \textit{Poletsky theory of discs in almost
complex manifolds}, (to appear in Complex Variables and Elliptic
Equations: An International Journal),
\bibitem[P]{p} N. Pali,  \textit{Fonctions plurisousharmoniques et courants positifs de type $(1, 1)$ sur une vari\'{e}t\'{e} presque
complexe},
Manuscripta Math. 118 (2005), no. 3, 311-337,

  \bibitem[Pl]{pl} S. Pli\`{s}, \textit{The Monge-Amp\`{e}re equation on almost complex manifolds}, ArXiv:1106.3356,
\bibitem[R]{r} R. Richberg, \textit{Stetige streng pseudokonvexe Funktionen}, Math. Annalen 175, 251-286
(1968),

\bibitem[S]{s}A. Sukhov,
\textit{Regularized maximum of strictly plurisubharmonic functions on an almost complex manifold}, ArXiv:1303.5312.
 
 \end{thebibliography}
\end{document}